\documentclass[12pt,reqno]{amsart}
\usepackage{amsmath, amsthm, amssymb}

\advance \topmargin by -\headheight
\advance \topmargin by -\headsep
     
\evensidemargin \oddsidemargin
\newtheorem{theorem}{Theorem}[section]

\newtheorem{corollary}[theorem]{Corollary}

\theoremstyle{definition}

\theoremstyle{remark}

\numberwithin{equation}{section}

\newcommand{\mmod}[1]{\,\,(\text{mod}\,\,#1)}

\def\calA{{\mathcal A}}

\def\calI{{\mathcal I}}

\def\calP{{\mathcal P}}

\def\calZ{{\mathcal Z}}

\def\dbN{{\mathbb N}}

\def\dbZ{{\mathbb Z}}\def\dbQ{{\mathbb Q}}

\def\grM{{\mathfrak M}}\def\grN{{\mathfrak N}}
\def\grS{{\mathfrak S}}\def\grP{{\mathfrak P}}

\def\grp{{\mathfrak p}}

\def\alp{{\alpha}} 
\def\bet{{\beta}}  
\def\gam{{\gamma}} \def\Gam{{\Gamma}}
\def\del{{\delta}}

\def\lam{{\lambda}}

\def\ome{{\omega}} 
\def\d{{\partial}}
\def\eps{\varepsilon}

\def\le{\leqslant} \def\ge{\geqslant}

\def\d{{\,{\rm d}}}

\begin{document}
\title[Waring's problem]{On Waring's problem:\\ two squares and three biquadrates}
\author[J. B. Friedlander]{John B. Friedlander$^*$}
\address{JBF: Department of Mathematics, University of Toronto, Toronto ON, 
M5S 2E4, Canada}
\email{frdlndr@math.toronto.edu}
\author[T. D. Wooley]{Trevor D. Wooley}
\address{TDW: School of Mathematics, University of Bristol, University Walk, Clifton, Bristol BS8 1TW, United Kingdom}
\email{matdw@bristol.ac.uk}
\thanks{$^*$Supported by NSERC grant A5123}
\date{}
\begin{abstract} We investigate sums of mixed powers involving two squares and three biquadrates. In particular, subject to the truth of the Generalised Riemann Hypothesis and the Elliott-Halberstam Conjecture, we show that all large natural numbers $n$ with $8\nmid n$, $n\not\equiv 2\mmod{3}$ and $n\not\equiv 14\mmod{16}$ are the sum of $2$ squares and $3$ biquadrates.\end{abstract}
\maketitle

\section{Introduction} Additive number theorists employed in the investigation of Waring's problem will recognise that the most challenging environment in which to ply their trade is that in which the sum of the reciprocals of the available exponents lies between $1$ and $2$. Indeed, while it is generally conjectured that, whenever $k_1^{-1}+\ldots +k_s^{-1}>1$, all large natural numbers $n$ satisfying appropriate congruence conditions should be representable in the form
\begin{equation}\label{1.1}
x_1^{k_1}+x_2^{k_2}+\ldots +x_s^{k_s}=n,
\end{equation}
in circumstances where $k_1^{-1}+\ldots +k_s^{-1}\le 2$ such has been established in only a handful of very special cases. Most of this work is restricted to problems containing a sum of three squares, or to problems containing two squares, two cubes, and various additional powers. Thus, on the one hand, Gauss \cite{Gau1801} and Hooley \cite{Hoo2000} respectively considered sums of three squares, and sums of three squares and a $k$th power. On the other hand, Linnik \cite{Lin1972} and Hooley \cite{Hoo1981b} investigated sums of two squares and three cubes. Very recently, Golubeva \cite{Gol1996, Gol2009} has shown that all large integers $n$ are represented as a sum of positive integral powers in the shape
\begin{equation}\label{1.2}
n=x_1^2+x_2^2+x_3^3+x_4^3+x_5^4+x_6^{16}+x_7^{4k+1}.
\end{equation}
Our goal in this paper is to investigate the more recalcitrant situation in which the cubes are replaced by biquadrates, the substantial escalation of difficulty requiring us to make use of two hypotheses that, although familiar to and widely believed by analytic number theorists, lie beyond the reach of current technology.\par

\begin{theorem}\label{theorem1.1}
Suppose that $n$ is a sufficiently large natural number with $8\nmid n$, $n\not\equiv 2\mmod{3}$ and $n\not\equiv 14\mmod{16}$.
Assume the Riemann Hypothesis for $L(s,\chi_{2n})$ and $L(s,\chi_{-2n})$ together with the Elliott-Halberstam Conjecture. Then $n$ is the sum of two squares and three biquadrates.
\end{theorem}

Here, we write $\chi_D(m)$ for the character $\left( \frac{D}{m}\right)$. Readers less familiar with the Elliott-Halberstam Conjecture will find the statement in (\ref{3.2}) below. Note that the sum of the reciprocals of the exponents in the representation problem
\begin{equation}\label{1.3}
n=x_1^2+x_2^2+x_3^4+x_4^4+x_5^4
\end{equation}
considered in Theorem \ref{theorem1.1} is $2-\frac{1}{4}$, whereas the corresponding sum in Golubeva's problem (\ref{1.2}) exceeds $2-\frac{1}{48}$. It is expected that all large natural numbers $n$ should be represented in the form (\ref{1.3}), the need to exclude congruence classes in the statement of Theorem \ref{theorem1.1} arising as an artefact of our methods. The reader will locate somewhat more relaxed congruential conditions in the penultimate paragraph of \S2 below. At the cost of augmenting the representation problem (\ref{1.3}) with an additional $k$th power, on the other hand, we are able to avoid the exclusion of any congruence class whatsoever.\par

\begin{corollary}\label{corollary1.2}
Let $k$ be a natural number. Then, under the same unproved assumptions as in the statement of Theorem \ref{theorem1.1}, every sufficiently large natural number $n$ is the sum of two squares, three biquadrates and a $k$th power.
\end{corollary}

The method that we use to establish Theorem \ref{theorem1.1} utilises representations of $n$ in the shape (\ref{1.3}) of special type. Although we obtain a lower bound of the expected magnitude for the number of these restricted representations, this falls short of the anticipated order of growth for their total number. If one is prepared to permit a small exceptional set of natural numbers $n$, then it is possible to not only dispense with the unproved assumptions, but even to obtain the anticipated asymptotic formula in the representation problem (\ref{1.3}). With this objective in mind, we introduce some notation with which to state our second theorem. Let $R_s(n)$ denote the number of representations of the positive number $n$ as the sum of two squares and $s$ biquadrates, and let
\begin{equation}\label{1.4}
\grS_s(n)=\sum_{q=1}^\infty \sum^q_{\substack{a=1\\ (a,q)=1}}\Bigl( q^{-1}\sum_{r=1}^qe(ar^2/q)\Bigr)^2\Bigl( q^{-1}\sum_{r=1}^qe(ar^4/q)\Bigr)^se(-na/q).
\end{equation}
Here, as usual, we write $e(z)$ for $e^{2\pi iz}$. We refer to a function $\psi(t)$ as being a {\it sedately increasing function} when $\psi(t)$ is a function of a positive variable $t$, increasing monotonically to infinity, and satisfying the condition that when $t$ is large, one has $\psi(t)=O(t^\del)$ for a positive number $\del$ sufficiently small in the ambient context. Finally, we write $E_s(X;\psi)$ for the number of integers $n$ with $1\le n\le X$ such that
\begin{equation}\label{1.5}
\left| R_s(n)-c_s\Gam(\tfrac{5}{4})^4\grS_s(n)n^{s/4}\right|>n^{s/4}\psi(n)^{-1},
\end{equation}
in which we write $c_3=\frac{2}{3}\sqrt{2}$ and $c_4=\frac{1}{4}\pi$. In \S4 we derive the relatively sharp upper bounds for $E_s(X;\psi)$ ($s=3,4$) presented in the following theorem.

\begin{theorem}\label{theorem1.2}
Suppose that $\psi(t)$ is a sedately increasing function. Then, for each $\eps>0$, one has
$$E_3(X;\psi)\ll X^{1/2+\eps}\psi(X)^2\quad \text{and}\quad E_4(X;\psi)\ll X^{1/4+\eps}\psi(X)^4,$$
where the implied constants may depend on $\eps$.
\end{theorem}

We conclude by noting that there is an extensive literature associated with problems of the shape (\ref{1.1}) involving mixed exponents of far wider generality than we have discussed herein. We refer the reader to Br\"udern \cite{Bru1987}, Hooley \cite{Hoo1981a} and Vaughan \cite{Vau1971} for a representative cross-section of such results. We remark also that Dietmann and Elsholtz \cite{DE2008} have recently shown that when $p$ is a prime number with $p\equiv 7\mmod{8}$, then $p^2$ cannot be written non-trivially in the form $p^2=x_1^2+x_2^2+x_3^4$, so that the exceptional set corresponding to two squares and a fourth power is necessarily rather large.\par

\section{Congruence conditions and a biquadratic identity}
Our method of proof of Theorem \ref{theorem1.1} rests on the application of the identity
\begin{equation}\label{2.1}
x^4+y^4+(x+y)^4=2(x^2+xy+y^2)^2.
\end{equation}
This has been extensively exploited elsewhere in the investigation of Waring's problem for sums of biquadrates (see for example work of Kawada and the second-named author \cite{KW1999}). Our basic strategy for proving Theorem \ref{theorem1.1} is to try to consider integers of the form $n-2p^2$, where $p$ runs through prime numbers congruent to $1$ modulo $3$. Of course, every such prime may be written in the form $x^2+xy+y^2$. Letting $\calA$ be the set of such integers with $p<\sqrt{n/2}$, we seek to sieve $\calA$ by the set of primes congruent to $3$ modulo $4$. If we are able to remove all such primes, we are left with a set of integers, 
each of which is the sum of two squares. In this situation, one obtains a representation of the shape $n=x_1^2+x_2^2+2p^2$, so that, in view of the identity (\ref{2.1}), one obtains the desired representation (\ref{1.3}) with $x_5=x_3+x_4$.\par

Examining the set $\calA$ it appears at first glance that, in order to reach integers which are the sum of two squares, one must sieve by all primes $\varpi\equiv 3\mmod{4}$ with $\varpi<n$. However, if one begins with a set of integers $n-2p^2$ known to be congruent to $1$ modulo $4$ then, since these integers have an even number of prime factors congruent to $3$ modulo $4$, it suffices to sieve up to $\sqrt{n}$. This makes an enormous difference. It is for this reason that we insist that $n$ lie in certain restricted congruence classes.\par

At this point we proceed in greater generality than is warranted for the proof of Theorem \ref{theorem1.1}. Our first step is to write $n$ uniquely 
in the form $n=16^hN_0$ where $16\nmid N_0$. We assume that
\begin{equation}\label{2.2}
N_0\not\equiv 14\mmod{16},\quad N_0\not\equiv 24\mmod{32}\quad \text{and}\quad N_0\not \equiv 2\mmod{3}.
\end{equation}
The argument splits into a small number of cases depending on the value of $N_0$. We define the integer $N$ and the set $\calA=\calA(N)$ as follows.\vskip.1cm

\noindent (i) When $N_0$ is congruent to $3$, $4$, $6$, $7$, $11$, $12$, or $15$ modulo $16$, we put $N=N_0$ and define
$$\calA(N)=\{ N-2p^2:\text{$p\equiv 1\mmod{3}$ and $p<\sqrt{N/2}$}\}.$$
\vskip.1cm
\noindent (ii) When $N_0$ is congruent to $1$, $2$, $5$, $9$, $10$, or $13$ modulo $16$, we put $N=N_0$ and define
$$\calA(N)=\{ N-32p^2:\text{$p\equiv 1\mmod{3}$ and $p<\sqrt{N/32}$}\}.$$
\vskip.1cm
\noindent (iii) When $N_0\equiv 8\mmod{16}$ we proceed as follows. Our assumption that $N_0\not\equiv 24\mmod{32}$ ensures that $N_0$ can be written in the shape $N_0=32N_1+8$ for some positive integer $N_1$. We put $N=8N_1+2$ and define
$$\calA(N)=\{ N-8p^2:\text{$p\equiv 1\mmod{3}$ and $p<\sqrt{N/8}$}\}.$$
\vskip.1cm

The set $\calA$ is defined in such a way that, in each case, it is indeed feasible that $\calA$ contain a sum of two squares. Indeed, in case (i) one finds that $N-2p^2$ is congruent to $1$ modulo $4$, or $2$ modulo $8$, or $4$ modulo $16$. Likewise, in case (ii) one may verify that $N-32p^2$ is congruent to $1$ modulo $4$ or $2$ modulo $8$, and in case (iii), instead $N-8p^2$ is congruent to $2$ modulo $8$. Observe next that, if we are successful in finding a sum of two squares $u^2+v^2$ in the set $\calA(N)$, then we are able to derive a representation of $n$ in the shape (\ref{1.3}). In order to confirm this, we recall that the implicit prime number $p$ may be written in the shape $p=a^2+ab+b^2$. Then in case (i) we have
$$n=(2^{2h}u)^2+(2^{2h}v)^2+(2^ha)^4+(2^hb)^4+(2^h(a+b))^4,$$
in case (ii) we have
$$n=(2^{2h}u)^2+(2^{2h}v)^2+(2^{h+1}a)^4+(2^{h+1}b)^4+(2^{h+1}(a+b))^4,$$
and in case (iii) we have
$$n=(2^{2h+1}u)^2+(2^{2h+1}v)^2+(2^{h+1}a)^4+(2^{h+1}b)^4+(2^{h+1}(a+b))^4.$$

\par We will see in \S3 that in each of the above cases, the set $\calA(N)$ does indeed contain a sum of two squares. In view of the above discussion, this suffices to complete the proof of Theorem \ref{theorem1.1}. Indeed, the above description implies the solubility of (\ref{1.3}) for any natural number $n$ such that, when $16^h\|n$ and $N_0=n/16^h$ is sufficiently large, then $N_0$ satisfies (\ref{2.2}). There are still however a small number of congruence classes for which one might expect a positive result but which remain untreated as an artefact of our method.\par

We complete this section by establishing Corollary \ref{corollary1.2}. Let $M$ be a large natural number. We note merely that, given an ordered pair $(r_1,r_2)\in \{0,1\}^2$, then for some integer $t\in\{1,16,33,48\}$, one has $t^k\equiv r_1\mmod{3}$ and $t^k\equiv r_2\mmod{16}$. For an appropriate choice of this integer $t$, one therefore finds that $n=M-t^k$ satisfies the hypotheses of Theorem \ref{theorem1.1}, and hence is the sum of two squares and three biquadrates. Thus $M$ is the sum of two squares, three biquadrates and a $k$th power.

\section{The proof of Theorem \ref{theorem1.1}} Once we are given the arguments of the previous section, the proof of Theorem \ref{theorem1.1} has been reduced to showing a positive lower bound for a certain sifting function. Fix the positive integer $N$. Define
$$S(\calA,\calP,z)=\sum_{\substack {m\in \calA\\ (m,P(z))=1}}1,$$
where $\calP$ denotes the set of primes congruent to $3$ modulo $4$ but not dividing $3N$, and
$$P(z)=\prod_{\substack{\varpi\in \calP\\ \varpi<z}}\varpi.$$
Then we seek to show that $S(\calA,\calP,z_0)$ is positive for the set $\calA=\calA(N)$ defined in the previous section in the various cases at hand, with $z_0=\sqrt{N}$. Our argument for doing this follows quite closely from that of the first-named author and Iwaniec in \cite[\S14.8]{FI2010}, so we shall simply sketch the main ideas, pointing out the very few points of difference.\par

As might be expected, the three subcases (i), (ii), (iii) described in \S2 also differ from each other very little. Indeed, since $\dbQ(\sqrt{2N})=\dbQ(\sqrt{32N})=\dbQ(\sqrt{8N})$ they all lead us to the same $L$-functions mentioned in the statement of the theorem. For the point of our discussion we shall specify that we are taking case (i).\par

We introduce first a lower bound semi-linear sieve $\{\lam_d:d<D,\, d| P(z)\}$ as described for example in \cite[Chapter 14]{FI2010}. We need to choose a ``level of distribution'' $D$ as large as possible yet such that the sieve remainder is small compared to the main term. The main term has a shape which depends mildly on $N$, but is certainly of size $\gg N^{1/2}/\log^2N$ (quadratic congruences modulo a prime have at most two solutions). The semi-linear sieve affords us a positive lower bound when the sieve parameter $s=(\log D)/(\log z)$ exceeds $1$. Hence, we can sieve by primes fairly close to $z_0$, say to $z=N^{1/2-\eps}$, provided we have available a level of distribution $D=N^{1/2-\eps/2}$.\par

Write $\rho(d)$ for the number of solutions of the congruence $\nu^2\equiv 2N\mmod{d}$. When $d|P(z)$, so that $d$ is odd and square-free, this is the same as the number of solutions of $2\nu^2\equiv N\mmod{d}$ and, moreover, one has $\rho(d)\le \tau(d)$. Note also that $(d,3)=1$ so that $\phi(3d)=2\phi(d)$ and, since we are running over primes congruent to $1$ modulo $3$, by the Chinese Remainder Theorem, the remainder in question is bounded by
\begin{equation}\label{3.1}
\sum_{\substack{d<D\\ d|P(z)}}\rho(d)\max_{(a,3d)=1}\left| \pi(\sqrt{N/2};3d,a)-\frac{\text{li}(\sqrt{N/2})}{\phi(3d)}\right|.
\end{equation}
The Elliott-Halberstam Conjecture states that for any $\eps>0$ one has the upper bound
\begin{equation}\label{3.2}
\sum_{q<Q}\max_{(b,q)=1}\left| \pi(x;q,b)-\frac{\text{li}(x)}{\phi(q)}\right|\ll x(\log x)^{-A}
\end{equation}
for any $A>0$, with $Q=x^{1-\eps}$, and an implied constant depending on $A$ and $\eps$. Take $Q=3D$.  The upper bound $\sqrt N (\log N)^{2-A/2}$, with $D=N^{1/2-\eps/2}$, for the sum in (\ref{3.1}) follows from this. In order to confirm this assertion, one has only to apply Cauchy's inequality, the Brun-Titchmarsh theorem, and the fact that 
$\sum_{d<D}\tau^2(d)/\phi(d)\ll \log^4 D$. 
We thus obtain the lower bound
\begin{equation}\label{3.3}
S(\calA,\calP,z)\gg \frac{\grS(N)\sqrt{N}}{(\log N)^{3/2}}\sqrt{s-1},
\end{equation}
where $s=(\log D)/(\log z)$ and where, in slight contrast to \cite[equation (14.78)]{FI2010}, 
$$\grS(N)=\left( \frac{L(1,\chi_{-2N})}{L(1,\chi_{2N})}\right)^{1/2}\prod_{\substack{\varpi|N\\ \varpi\equiv 3\mmod{4}}}(1+1/\varpi).$$
The Riemann Hypothesis for these two $L$-functions is used in this 
evaluation.\par

Since $z$ is somewhat smaller than $\sqrt{N}$ our sifting function $S(\calA,\calP,z)$ may count some integers having two prime factors $\varpi_1\equiv \varpi_2\equiv 3\mmod{4}$ with $z<\varpi_1\le \varpi_2\le \sqrt{N}$. We wish to subtract the contribution $T(\calA,\calP,z)$ coming from those integers and require for this an upper bound which is smaller than the lower bound in (\ref{3.3}). An upper bound is provided for this in \cite{FI2010}, but only with the aid of deep results \cite{DFI2012} on Weyl sums for quadratic roots. The problem of bounding $T$ involves the replacement of each of the variables $p$, $\varpi_1$, $\varpi_2$ by an upper-bound sieve while carefully maintaining control of the close proximity of $\varpi_1$ and $\varpi_2$. After implementation of the sieves the issue reduces to a counting problem for integer points on a certain hyperboloid. The difference in our case, aside from slightly different coefficients in the equation of the hyperboloid, comes from the condition $p\equiv 1\mmod{3}$ in our set $\calA$. Since we need only an upper bound for $T$ we can simply throw this condition away and appeal to the argument for \cite[equation (14.98)]{FI2010}. This yields the upper bound
$$T(\calA,\calP,z)\ll \frac{\grS(N)\sqrt{N}}{(\log N)^{3/2}}(s-1)^{3/2},$$
which, combined with (\ref{3.3}) and choosing $s$ sufficiently close to $1$, yields
$$S(\calA,\calP,\sqrt{N})=S(\calA,\calP,z)-T(\calA,\calP,z)\gg \grS(N)\sqrt{N}(\log N)^{-3/2}.$$
This completes the proof of Theorem \ref{theorem1.1}. 

\par We note that, rather than restricting in $\calA$ to primes congruent to $1$ modulo $3$, one might try to employ the larger set
$$\calA = \{N-2(y_1^2+y_1y_2+y_2^2)^2: y_1^2+y_1y_2+y_2^2\le \sqrt{N/2}\}.$$
This would require an Elliott-Halberstam like assumption not for primes but for norms in the quadratic field $\dbQ(\sqrt{-3})$. Although not so well known and no doubt very difficult, this assumption for a divisor-like function could conceivably be less distant than its better known counterpart for primes.

\section{Exceptional sets} Our approach to the job of proving Theorem \ref{theorem1.2} is motivated by work of the second-named author joint with Kawada \cite{KW2010}. Before we launch our proof in earnest, a word is in order on the convention adopted in this section concerning the use of the number $\eps$. Whenever $\eps$ appears in a statement, either implicitly or explicitly, we assert that the statement holds for each $\eps>0$. Note that the ``value'' of $\eps$ may consequently change from statement to statement.\par

We take $s$ to be either $3$ or $4$. Suppose that $X$ is a large positive number, and let $\psi(t)$ be a sedately increasing function. We denote by $\calZ_s(X)$ the set of integers $n$ with $X/2<n\le X$ for which the lower bound (\ref{1.5}) holds, and we abbreviate $\text{card}(\calZ_s(X))$ to $Z_s$. Write $P_k$ for $[X^{1/k}]$, and define the exponential sum $f_k(\alp)$ by
$$f_k(\alp)=\sum_{1\le x\le P_k}e(\alp x^k).$$
Also, when $Q$ is a positive number, let $\grM(Q)$ denote the union of the intervals
$$\grM(q,a)=\{ \alp\in [0,1):|q\alp-a|\le QX^{-1}\},$$
with $0\le a\le q\le Q$ and $(a,q)=1$. Then, when $1\le Q\le X^{1/2}$, we put $\grN(Q)=\grM(2Q)\setminus \grM(Q)$. In order to ensure that each $\alp\in [0,1)$ is associated with a uniquely defined arc $\grM(q,a)$, we adopt the convention that when $\alp$ lies in more than one arc $\grM(q,a)\subseteq \grM(Q)$, then it is declared to lie in the arc for which $q$ is least. Finally, let $\nu$ be a sufficiently small positive number, and write $W=X^\nu$. We then take $\grP$ to be the union of the intervals
$$\grP(q,a)=\{\alp \in [0,1):|\alp -a/q|\le WX^{-1}\},$$
with $0\le a\le q\le W$ and $(a,q)=1$. We set $\grp=[0,1)\setminus \grP$
and note that, by Dirichlet's theorem, the minor arcs $\grp$ can be covered by a dyadic union of the arcs $\grN(Q)$ with $X^{\nu} \le Q \le X^{1/2}$. We remark that the arcs $\grP(q,a)$ are wider than the corresponding arcs $\grM(q,a)$ comprising $\grM(W)$, but with a width which is independent of $q$. The latter property eases technical complications associated with the asymptotic analysis of the major arc contribution.\par

We begin by sketching how the methods of \cite[Chapter 4]{Vau1997} lead to the asymptotic formula
\begin{equation}\label{4.1}
\int_\grP f_2(\alp)^2f_4(\alp)^se(-n\alp)\d\alp =\frac{\Gam(\tfrac{3}{2})^2\Gam(\tfrac{5}{4})^s}{\Gam(\tfrac{s}{4}+1)}\grS_s(n)n^{s/4}+O(n^{s/4-\tau}),
\end{equation}
for a suitably small positive number $\tau$. Here, the singular series $\grS_s(n)$ is that defined in (\ref{1.4}). We note in this context that when $s$ is $3$ or $4$, one may verify by exploiting standard properties of the $\Gam$-function that the constant $c_s$ introduced following (\ref{1.5}) satisfies
$$c_s\Gam(\tfrac{5}{4})^4=\frac{\Gam(\tfrac{3}{2})^2\Gam(\tfrac{5}{4})^s}{\Gam(\tfrac{s}{4}+1)}.$$

\par Write
$$S_k(q,a)=\sum_{r=1}^qe(ar^k/q)\quad \text{and}\quad v_k(\bet)=\int_0^{P_k}e(\bet \gam^k)\d\gam ,$$
and for the moment define $f_k^*(\alp)$ for $\alp\in \grP(q,a)\subseteq \grP$ by putting
\begin{equation}\label{4.2}
f_k^*(\alp)=q^{-1}S_k(q,a)v_k(\alp-a/q).
\end{equation}
Then it follows from \cite[Theorem 4.1]{Vau1997} that whenever $\alp \in \grP(q,a)\subseteq \grP$, one has
$$f_k(\alp)-f_k^*(\alp)\ll q^{1/2+\eps}\ll W^{1/2+\eps}.$$
The measure of $\grP$ is $O(W^3X^{-1})$, and so it follows that
$$\int_\grP f_2(\alp)^2f_4(\alp)^se(-n\alp)\d\alp -\int_\grP f_2^*(\alp)^2f_4^*(\alp)^se(-n\alp)\d\alp \ll W^4X^{(s-1)/4}.$$
But a routine computation confirms that
$$\int_\grP f_2^*(\alp)^2f_4^*(\alp)^se(-n\alp)\d\alp =\grS_s(n;W)J_s(n;W),$$
where
$$\grS_s(n;W)=\sum_{1\le q\le W}A_s(q;n),$$
in which we have written
\begin{equation}\label{4.4}
A_s(q;n)=\sum^q_{\substack{a=1\\ (a,q)=1}}(q^{-1}S_2(q,a))^2(q^{-1}S_4(q,a))^se(-na/q),
\end{equation}
and
$$J_s(n;W)=\int_{-WX^{-1}}^{WX^{-1}}v_2(\bet)^2v_4(\bet)^se(-n\bet)\d\bet .$$
Thus we may conclude at this stage that
\begin{equation}\label{4.5}
\int_\grP f_2(\alp)^2f_4(\alp)^se(-n\alp)\d\alp=\grS_s(n;W)J_s(n;W)+O(X^{s/4-\tau}).
\end{equation}

\par The estimate
$$v_k(\bet)\ll P_k(1+|\bet|X)^{-1/k},$$
available from \cite[Theorem 7.3]{Vau1997}, ensures that the singular integral $J_s(n;W)$ converges absolutely as $W\rightarrow \infty$, and the discussion concluding \cite[Chapter 4]{Dav1963} readily leads to the relation
\begin{equation}\label{4.6}
J_s(n;W)=\frac{\Gam(\tfrac{3}{2})^2\Gam(\tfrac{5}{4})^s}{\Gam(\tfrac{s}{4}+1)}n^{s/4}+O(n^{s/4-\tau}).
\end{equation}

\par In order to analyse the corresponding truncated singular series $\grS_s(n;W)$, we begin by defining the multiplicative function $w_k(q)$ by taking
$$w_k(p^{uk+v})=\begin{cases} kp^{-u-1/2},&\text{when $u\ge 0$ and $v=1$,}\\
p^{-u-1},&\text{when $u\ge 0$ and $2\le v\le k$.}\end{cases}$$
Then according to \cite[Lemma 3]{Vau1986}, whenever $a\in \dbZ$ and $q\in \dbN$ satisfy $(a,q)=1$, one has $q^{-1}S_k(q,a)\ll w_k(q)$. We therefore deduce from (\ref{4.4}) that when $q=p^h$, then
$$q^{1/4}A_s(q;n)\ll q^{5/4}w_2(q)^2w_4(q)^s\ll p^{-\max\{5/4,h/2\}},$$
whence
$$\sum_{h=1}^\infty p^{h/4}|A_s(p^h;n)|\ll p^{-5/4}+\sum_{h=3}^\infty p^{-h/2}\ll p^{-5/4}.$$
By multiplicativity of the function $A_s(q;n)$, therefore, one finds that for a suitable positive number $B$,
\begin{align*}
\sum_{q=1}^\infty (q/W)^{1/4}|A_s(q;n)|&=W^{-1/4}\prod_p\Bigl( 1+\sum_{h=1}^\infty p^{h/4}|A_s(p^h;n)|\Bigr)\\
&\ll W^{-1/4}\prod_p (1+Bp^{-5/4})\ll W^{-1/4}.
\end{align*}
Thus we deduce that
$$\sum_{q>W}|A_s(q;n)|\ll W^{-1/4},$$
hence that the singular series $\grS_s(n)$ defined in (\ref{1.4}) converges absolutely, and in particular that
$$\grS_s(n)-\grS_s(n;W)\ll W^{-1/4}.$$
We may therefore conclude from (\ref{4.5}) and (\ref{4.6}) that the relation (\ref{4.1}) does indeed hold.\par

We may now proceed to examine the exceptional set itself. For each integer $n\in \calZ_s(X)$, it follows from (\ref{1.5}) via orthogonality that
$$\left| \int_0^1f_2(\alp)^2f_4(\alp)^se(-n\alp)\d\alp -c_s\Gam(\tfrac{5}{4})^4\grS_s(n)n^{s/4}\right| >\tfrac{1}{2}X^{s/4}\psi(X)^{-1}.$$
In view of the relation (\ref{4.1}), therefore, we obtain the lower bound
$$\left| \int_\grp f_2(\alp)^2f_4(\alp)^se(-n\alp)\d\alp \right| >\tfrac{1}{4}X^{s/4}\psi(X)^{-1},$$
whence
$$\sum_{n\in \calZ_s(X)}\left| \int_\grp f_2(\alp)^2f_4(\alp)^se(-n\alp)\d\alp \right| \gg Z_sX^{s/4}\psi(X)^{-1}.$$
There exist complex numbers $\eta_n=\eta_n(s)$, with $|\eta_n|=1$, satisfying the condition that for each $n\in \calZ_s(X)$, one has
$$\Bigl|\int_\grp f_2(\alp)^2f_4(\alp)^se(-n\alp)\d\alp \Bigr|=\eta_n(s)\int_\grp f_2(\alp)^2f_4(\alp)^se(-n\alp)\d\alp .$$
Consequently, with the exponential sum $K_s(\alp)$ defined by
$$K_s(\alp)=\sum_{n\in \calZ_s(X)}\eta_n(s)e(-n\alp),$$
one finds that
\begin{align}
Z_sX^{s/4}\psi(X)^{-1}&\ll \sum_{n\in \calZ_s(X)}\eta_n(s)\int_\grp f_2(\alp)^2f_4(\alp)^se(-n\alp)\d\alp \notag \\
&\le \int_\grp |f_2(\alp)^2f_4(\alp)^sK_s(\alp)|\d\alp .\label{4.7}
\end{align}

We next derive an upper bound for the integral on the right hand side of (\ref{4.7}). Let $Q$ be a positive real number with $W\le Q\le X^{1/2}$. Write
\begin{equation}\label{4.8}
\calI=\int_0^1|f_2(\alp)^2f_4(\alp)^4|\d\alp,
\end{equation}
and recall the familiar estimate
\begin{equation}\label{4.9}
\calI\ll X^{1+\eps}
\end{equation}
available, for example, from \cite[Exercise 6 of Chapter 2]{Vau1997}. Define
$$\ome_s=\begin{cases} 0,&\text{when $s=3$,}\\
1,&\text{when $s=4$,}\end{cases}$$
and put
\begin{equation}\label{4.10}
I_s=\int_0^1|f_4(\alp)^{2s-6}K_s(\alp)^2|\d\alp .
\end{equation}
Then from Parseval's identity and \cite[Lemma 2.1]{KW2010}, one obtains the bound
\begin{equation}\label{4.11}
I_s\ll P_4^{s-3}Z_s+\ome_sP_4^{1/2+\eps}Z_s^{3/2}\quad (s=3,4).
\end{equation}

\par We modify our notation for major arc approximants in order to facilitate the subsequent analysis. Define now $f_k^*(\alp)$ via (4.2) for $\alp \in \grM(q,a)\subseteq \grM(Q)$. Then it follows from \cite[Theorem 4.1]{Vau1997} that whenever $\alp \in \grN(Q)$, then
$$f_k(\alp)-f_k^*(\alp)\ll Q^{1/2+\eps},$$
and from \cite[Theorems 4.2 and 7.3]{Vau1997} that
$$\sup_{\alp \in \grN(Q)}|f_k^*(\alp)|\ll P_kQ^{-1/k}.$$
Thus, in particular, we find that under our present hypotheses concerning $Q$, whenever $\alp \in \grN(Q)$ one has
\begin{align*}
f_2(\alp)f_4(\alp)-f_2^*(\alp)f_4^*(\alp)&\ll Q^{1/2+\eps}(|f_2^*(\alp)|+|f_4^*(\alp)|)+Q^{1+\eps}\\
&\ll X^\eps (P_2+P_4Q^{1/4}+X^{1/2})\ll X^{1/2+\eps}.
\end{align*}
On writing
\begin{equation}\label{4.12}
U_s=\int_{\grN(Q)}|f_2(\alp)f_2^*(\alp)f_4(\alp)^{s-1}f_4^*(\alp)K_s(\alp)|\d\alp ,
\end{equation}
we therefore deduce that
\begin{equation}\label{4.13}
\int_{\grN(Q)}|f_2(\alp)^2f_4(\alp)^sK_s(\alp)|\d\alp\ll U_s+X^{1/2+\eps} \int_0^1|f_2(\alp)f_4(\alp)^{s-1}K_s(\alp)|\d\alp .
\end{equation}

\par On recalling (\ref{4.8}) and (\ref{4.10}), and applying Schwarz's inequality, we discern that
$$\int_0^1|f_2(\alp)f_4(\alp)^{s-1}K_s(\alp)|\d\alp \le I_s^{1/2}\calI^{1/2}.$$
Then from (\ref{4.9}) and (\ref{4.11}), we see that when $s$ is $3$ or $4$ one has
\begin{equation}\label{4.14}
\int_0^1|f_2(\alp)f_4(\alp)^{s-1}K_s(\alp)|\d\alp \ll X^{1/2+\eps}(P_4^{s-3}Z_s+\ome_sP_4^{1/2}Z_s^{3/2})^{1/2}.
\end{equation}

Turning our attention next to the integral $U_s$ defined in (\ref{4.12}), we begin by observing that \cite[Theorems 4.2 and 7.3]{Vau1997} yield the bound
$$f_2^*(\alp)f_4^*(\alp)\ll X^{3/4}G(\alp)^{3/4},$$
where we write $G(\alp)$ for the function defined by taking
$$G(\alp)=(q+X|q\alp-a|)^{-1},$$
when $\alp\in \grM(q,a)\subseteq \grM(Q)$, and $0$ otherwise. Then on recalling (\ref{4.8}), an application of Schwarz's inequality reveals that
\begin{equation}\label{4.15}
U_s\ll X^{3/4}\Bigl( \sup_{\alp \in \grN(Q)}|G(\alp)|\Bigr)^{1/4}\calI^{1/2}V_s^{1/2},
\end{equation}
where
$$V_s=\int_{\grN(Q)}G(\alp)\Psi(\alp)\d\alp,$$
in which we put $\Psi(\alp)=|f_4(\alp)^{s-3}K_s(\alp)|^2$.\par

Writing
$$\Psi(\alp)=\sum_{|h|\le 2X}\psi_he(\alp h),$$
we find from \cite[Lemma 2]{Bru1988} that
$$V_s\ll X^{\eps-1}\Bigl(Q\psi_0+\sum_{h\ne 0}|\psi_h|\Bigr).$$
On putting
$$K_s^*(\alp)=\sum_{n\in \calZ_s(X)}e(-n\alp),$$
we see from orthogonality that
$$|\psi_h|\le \int_0^1|f_4(\alp)^{s-3}K_s^*(\alp)|^2e(-\alp h)\d\alp,$$
and hence
$$V_s\ll X^{\eps-1}\Bigl( Q\int_0^1|f_4(\alp)^{2s-6}K_s(\alp)^2|\d\alp +f_4(0)^{2s-6}K_s^*(0)^2\Bigr) .$$
Then one sees by means of (\ref{4.10}) and (\ref{4.11}) that
$$V_s\ll X^{\eps-1}\Bigl( Q(P_4^{s-3}Z_s+\ome_sP_4^{1/2}Z_s^{3/2})+P_4^{2s-6}Z_s^2\Bigr) .$$
Then in view of (\ref{4.9}) and the trivial estimate
$$\sup_{\alp \in \grN(Q)}|G(\alp)|^{1/4}\ll Q^{-1/4},$$
we deduce from (\ref{4.15}) that whenever $X^\nu\le Q\le X^{1/2}$, one has
\begin{align*}
U_s&\ll X^{3/4+\eps}Q^{-1/4}\left( Q(P_4^{s-3}Z_s+\ome_sP_4^{1/2}Z_s^{3/2})+P_4^{2s-6}Z_s^2\right) ^{1/2}\\
&\ll X^{s/4+\eps}\left( X^{1/8}(X^{(3-s)/8}Z_s^{1/2}+\ome_sX^{(5-2s)/16}Z_s^{3/4})+X^{-\nu/5}Z_s\right).
\end{align*}
Finally, on recalling (\ref{4.13}) and (\ref{4.14}), and keeping in mind that $s$ is either $3$ or $4$, we reach the upper bound
\begin{align*}\int_{\grN(Q)}|f_2(\alp)^2&f_4(\alp)^sK_s(\alp)|\d\alp \\
&\ll X^{s/4+\eps}(X^{(5-s)/8}Z_s^{1/2}+\ome_sX^{(9-2s)/16}Z_s^{3/4}+X^{-\nu/5}Z_s).
\end{align*}

\par By covering the minor arcs $\grp$ by a dyadic union of arcs $\grN(Q)$, therefore, we discover from (\ref{4.7}) that
$$Z_sX^{s/4}\psi(X)^{-1}\ll X^{s/4+\eps}(X^{(5-s)/8}Z_s^{1/2}+\ome_sX^{(9-2s)/16}Z_s^{3/4}+X^{-\nu/5}Z_s).$$
Thus, provided that $\psi(X)\ll X^\del$ for a sufficiently small positive number $\del$, as we are permitted to suppose, we disentangle the estimate
$$Z_s\ll X^{(5-s)/4+\eps}\psi(X)^2+\ome_sX^{(9-2s)/4+\eps}\psi(X)^4,$$
which in turn yields the bounds
$$Z_3\ll X^{1/2+\eps}\psi(X)^2\quad \text{and}\quad Z_4\ll X^{1/4+\eps}\psi(X)^4.$$
The conclusion of Theorem \ref{theorem1.2} follows by summing the contributions from $\calZ_s(X)$ over dyadic intervals in $X$ so as to cover the interval $[1,X]$.

\bibliographystyle{amsbracket}

\begin{thebibliography}{18}

\bibitem{Bru1987}
J. Br\"udern, \emph{Sums of squares and higher powers}, J. London Math. Soc. (2) \textbf{35} (1987), 233--243.

\bibitem{Bru1988}
J. Br\"udern, \emph{A problem in additive number theory}, Math. Proc. Cambridge Philos. Soc. \textbf{103} (1988), 27--33.

\bibitem{Dav1963}
H. Davenport, \emph{Analytic methods for Diophantine equations and Diophantine inequalities}, Second edition, Cambridge University Press, Cambridge, 2005.

\bibitem{DE2008}
R. Dietmann and C. Elsholtz, \emph{Sums of two squares and one biquadrate}, Funct. Approx. Comment. Math. \textbf{38} (2008), 233--234.

\bibitem{DFI2012}
W. Duke, J. B. Friedlander and H. Iwaniec, \emph{Weyl sums for quadratic roots}, Internat. Math. Res. Notices (2012), no. 11, 2493--2549; \emph{Corrigendum}, Internat. Math. Res. Notices (2012), no. 11, 2646-2648.

\bibitem{FI2010}
J. B. Friedlander and H. Iwaniec, \emph{Opera de Cribro}, American Mathematical Society, Providence, RI, 2010.

\bibitem{Gau1801}
C. F. Gauss, \emph{Disquisitiones arithmeticae}, Leipzig, 1801.

\bibitem{Gol1996}
E. P. Golubeva, \emph{On nonhomogeneous Waring equations}, Zap. Nauchn. Sem. S.-Peterburg. Otdel. Mat. Inst. Steklov. (POMI) \textbf{226} (1996), Anal. Teor. Chisel i Teor. Funktsii. \textbf{13}, 65--68, 236; translation in
J. Math. Sci. (New York) \textbf{89} (1998), 955--957. 

\bibitem{Gol2009}
E. P. Golubeva, \emph{A bound for the representability of large numbers by ternary forms, and nonhomogeneous Waring equations}, Zap. Nauchn. Sem. S.-Peterburg. Otdel. Mat. Inst. Steklov. (POMI) \textbf{357} (2008), Analiticheskaya Teoriya Chisel i Teoriya Funktsii. \textbf{23}, 5--21, 224; translation in J. Math. Sci. (New York) \textbf{157} (2009), 543--552.

\bibitem{Hoo1981a}
C. Hooley, {\it On a new approach to various problems of Waring's type}, Recent progress in analytic number theory, Vol. 1 (Durham, 1979), pp. 127--191, Academic Press, London-New York, 1981. 

\bibitem{Hoo1981b}
C. Hooley, \emph{On Waring's problem for two squares and three cubes}, J. Reine Angew. Math. \textbf{328} (1981), 161--207.

\bibitem{Hoo2000}
C. Hooley, \emph{On Waring's problem for three squares and an $l$th power}, Asian J. Math. \textbf{4} (2000), 885--903.

\bibitem{KW1999}
K. Kawada and T. D. Wooley, \emph{Sums of fourth powers and related topics}, J. Reine Angew. Math. \textbf{512} (1999), 173--223.

\bibitem{KW2010}
K. Kawada and T. D. Wooley, \emph{Davenport's method and slim exceptional sets: the asymptotic formulae in Waring's problem}, Mathematika \textbf{56} (2010), 305--321.

\bibitem{Lin1972}
Yu. V. Linnik, \emph{Additive problems involving squares, cubes and almost primes}, Acta Arith. \textbf{21} (1972), 413--422. 

\bibitem{Vau1971}
R. C. Vaughan, \emph{On sums of mixed powers}, J. London Math. Soc. (2) \textbf{3} (1971), 677--688. 

\bibitem{Vau1986}
R. C. Vaughan, \emph{On Waring's problem for smaller exponents}, Proc. London Math. Soc. (3) \textbf{52} (1986), 445--463.

\bibitem{Vau1997}
R. C. Vaughan, \emph{The Hardy-Littlewood method}, 2nd edition, Cambridge University Press, Cambridge, 1997.

\end{thebibliography}
\providecommand{\bysame}{\leavevmode\hbox to3em{\hrulefill}\thinspace}

\end{document}